\documentclass{amsart}
\usepackage{amsmath}
\usepackage{amsfonts}
\usepackage{graphics}
\usepackage{epsfig}
\usepackage{amssymb}
\usepackage{amscd}

\newtheorem{thm}{Theorem}[section]

\newtheorem{lem}[thm]{Lemma}

\theoremstyle{definition}

\newtheorem*{ack}{Acknowledgement}
\theoremstyle{remark}
\newtheorem{rem}[thm]{Remark}

\numberwithin{equation}{section}
\numberwithin{figure}{section}

\def\Hom{{\text{\rm{Hom}}}}
\def\End{{\text{\rm{End}}}}
\def\trace{{\text{\rm{trace}}}}
\def\tr{{\text{\rm{tr}}}}
\def\rchi{{\hbox{\raise1.5pt\hbox{$\chi$}}}}
\def\Aut{{\text{\rm{Aut}}}}

\def\isom{\cong}
\def\tensor{\otimes}
\def\dsum{\oplus}

\def\reg{{\text{\rm{reg}}}}

\begin{document}
\large
\setcounter{section}{-1}

\title[Homomorphisms from 
a surface group into a finite group]{A 
generating function of 
the number of homomorphisms from a 
surface group into a finite group$^\dagger$}
\author[Motohico Mulase]{Motohico Mulase$^1$}  
\address{
Department of Mathematics\\
University of California\\
Davis, CA 95616--8633}
\email{mulase@math.ucdavis.edu}
\author[Josephine Yu]{Josephine T.\ Yu$^2$}  
\address{
Department of Mathematics\\
University of California\\
Davis, CA 95616--8633}
\email{yujt@math.ucdavis.edu}
\begin{abstract} A generating function of 
the number of homomorphisms
from the fundamental group of a compact 
oriented or
non-orientable surface without 
boundary into a
finite group is obtained
in terms of an integral over
a real group algebra. We calculate 
the number of homomorphisms 
using the decomposition of the group algebra into
irreducible factors. This gives
a new proof of the classical formulas
of Frobenius, Schur, and Mednykh.
\end{abstract}

\date{November 5, 2002}
\thanks{$^\dagger$This article is based on the first
author's talks at CIMAT, Guanajuato, and
MSRI, Berkeley, in September, 2002,
and serves as an announcement
of \cite{MY}.}
\thanks{$^1$Research supported   by 
NSF grant DMS-9971371
and UC Davis.}
\thanks{$^2$Research supported   by NSF grant VIGRE 
DMS-0135345 and UC Davis.}
\maketitle

\allowdisplaybreaks

%\tableofcontents

\section{Introduction}

\noindent
Let $S$ be a compact oriented or non-orientable surface
without boundary, and $\rchi(S)$ its Euler characteristic. 
The subject of our study is a generating function of the number
$|\Hom(\pi_1(S),G)|$ of homomorphisms from the 
fundamental group of $S$ into a finite group $G$. 
We give a generating function  in terms
of a \emph{non-commutative  integral}
Eqn.(\ref{eq:CGasymptotic}) or  Eqn.(\ref{eq:realgr}),
according to the orientability of $S$.
 The idea of such
integrals comes from random matrix theory.
Our integrals can be thought of as
 a generalization  of real symmetric, complex
hermitian, and quaternionic self-adjoint 
matrix integrals.

The graphical expansion methods for
real symmetric \cite{BIPZ, GHJ}
and complex hermitian \cite{BIZ}
 matrix integrals are generalized in \cite{MW}
for quaternionic self-adjoint matrix integrals. 
The technique developed in \cite{MW} is further
generalized in \cite{MY} to 
the integrals over matrices
with values in non-commutative $*$-algebras.
In this article we consider 
$1\times 1$ matrix integrals
over group algebras. Surprisingly, the graphical 
expansion of the integral gives a generating function
of  $|\Hom(\pi_1(S),G)|$ for all closed surfaces.

\section{Counting Formulas}

\noindent
Computation of our generating functions 
Eqn.(\ref{eq:CGasymptotic}) and Eqn.(\ref{eq:realgr})
yields a new proof of  the 
following classical counting formulas:

\begin{thm}[Mednykh \cite{Med}]
\label{thm:orientable}
Let $G$ be a finite group of order $|G|$, and $\hat{G}$
the set of all complex irreducible representations of $G$. 
By $V_\lambda$ we denote the irreducible representation 
parameterized by $\lambda\in \hat{G}$. Then for every
compact Riemann surface $S$, we have

\begin{equation}
\label{eq:orihomcount}
\sum_{\lambda\in\hat{G}} (\dim V_\lambda)^{\rchi(S)}
= |G|^{\rchi(S)-1}\cdot |\Hom(\pi_1(S),G)|\;,
\end{equation}
where $\rchi(S)$ is the Euler characteristic of the surface
$S$.
\end{thm}

\begin{rem}
{F}or a surface of genus $0$, 
Eqn.(\ref{eq:orihomcount}) reduces to the classical
formula
$$
\sum_{\lambda\in\hat{G}} (\dim V_\lambda) ^2 =|G|\;.
$$
 Since $\dim V_\lambda$ is a
divisor of  $|G|$, we note that the expression
$$
|G|^{-\rchi(S)} \sum_{\lambda\in\hat{G}}
 (\dim V_\lambda) ^{\rchi(S)} = 
|\Hom(\pi_1(S),G)| \big/|G|
$$
is an integer for a surface of positive genus. 
\end{rem}

{F}or a non-orientable surface, there is another formula:

\begin{thm}[Frobenius-Schur \cite{FS}]
\label{thm:nonorihomcount}
Let us
decompose $\hat{G}$ into three disjoint subsets 
according to the Frobenius-Schur indicator: 
\begin{equation}
\label{eq:3reps}
\begin{split}
\hat{G}_1 &= \bigg\{\lambda\in\hat{G} 
\;\bigg|\;\frac{1}{|G|}\sum_{\gamma\in G} 
\rchi_\lambda (\gamma^2) = 1 
\bigg\}\\
\hat{G}_2 &= \bigg\{\lambda\in\hat{G} 
\;\bigg|\;\frac{1}{|G|}\sum_{\gamma\in G} 
\rchi_\lambda (\gamma^2) = 0 
\bigg\}\\
\hat{G}_3 &= \bigg\{\lambda\in\hat{G} 
\;\bigg|\;\frac{1}{|G|}\sum_{\gamma\in G} 
\rchi_\lambda (\gamma^2) = -1 
\bigg\} \; .
\end{split} 
\end{equation}
Then we have
\begin{equation}
\label{eq:nonorihomcount}
\sum_{\lambda\in\hat{G}_1} 
(\dim_\mathbb{C} V_\lambda) ^{\rchi(S)} + 
\sum_{\lambda\in\hat{G}_3} 
(-\dim_\mathbb{C} V_\lambda) ^{\rchi(S)}
= |G|^{\rchi(S)-1}\cdot  |\Hom(\pi_1(S),G)| \;,
\end{equation}
where $S$ is now an arbitrary compact 
non-orientable surface
without boundary.
\end{thm}

\begin{rem}
If we take $S=\mathbb{R}P^2$, then 
$\rchi(\mathbb{R}P^2)=1$ and $\pi_1(\mathbb{R}P^2)
= \mathbb{Z}/2\mathbb{Z}$, and the formula 
reduces to a well-known formula \cite{Isaacs, Serre}
$$
\sum_{\lambda\in\hat{G}_1} \dim_\mathbb{C} V_\lambda
- \sum_{\lambda\in\hat{G}_3} \dim_\mathbb{C} V_\lambda
= \text{ the number of involutions in } G.
$$
In particular, if every complex irreducible 
representation of
$G$ is defined over $\mathbb{R}$, then 
$\hat{G} = \hat{G}_1$ and the same formula
(\ref{eq:orihomcount}) holds for an arbitrary 
closed surface $S$, 
orientable or 
non-orientable.
This class of groups includes symmetric groups
$\mathfrak{S}_n$. 
\end{rem}

\section{Orientable Case}

\noindent
The generating function of $|\Hom(\pi_1(S),G)|$
for an oriented surface $S$ comes from analysis
of hermitian matrix integrals. The prototype of the 
asymptotic expansion formula is

\begin{equation}
\label{eq:proto}
\log  \int_{\mathcal{H}_{N,\mathbb{C}}}
e^{-\frac{1}{2}N\tr X^2} 
e^{N\sum_{j\ge 3} \frac{t_j}{j}\tr X^j}d\mu(X)
=
\sum_{\substack{\Gamma \text{ connected}\\
\text{Ribbon graph}}}
\frac{1}{|\Aut_{R} \Gamma|} N^{\rchi(S_\Gamma)} 
\prod_{j\ge 3} t_j ^{v_j (\Gamma)}
\end{equation}
for an $N\times N$ hermitian matrix integral
\cite{BIZ}, where 
$\mathcal{H}_{N,\mathbb{C}}$ denotes the 
space of hermitian matrices of size $N$, and
$e^{-\frac{N}{2}\tr X^2} d\mu(X)$  the normalized
probability measure on 
$\mathcal{H}_{N,\mathbb{C}} = \mathbb{R}^{N^2}$.
A \emph{ribbon graph} is a graph with a cyclic order assigned
to each vertex. Equivalently, it is a graph $\Gamma$ that is 
drawn on a closed oriented surface $S$ such that the complement
$S\setminus \Gamma$ is the disjoint union of open disks.
We use $f(\Gamma)$ to denote the number of these disks 
(or \emph{faces}), and
$v_j (\Gamma)$ the number of $j$-valent vertices of $\Gamma$.
Let 
$$
v(\Gamma) = \sum_j v_j(\Gamma) 
\qquad \text{and} \qquad
e(\Gamma) = \frac{1}{2}\; \sum_j j v_j(\Gamma)
$$ 
be the number of vertices
and edges of $\Gamma$, respectively. Then the genus 
$g(S)$ of
the surface $S$ is given by a  formula for 
Euler characteristic
$$
\rchi(S) = 2-2g(S) = v(\Gamma) - e(\Gamma) + f(\Gamma).
$$
The automorphism group $\Aut_R(\Gamma)$ of a ribbon graph
$\Gamma$ consists of automorphisms of  
the cell-decomposition of $S$ that is determined by the graph.
We refer to \cite{MP1998} for precise definition
of ribbon graphs and their automorphism groups.

The asymptotic formula is an equality in the ring
$$
\mathbb{Q}(N)[[t_3,t_4,t_5,\dots]]
$$
of formal power series in an infinite number of variables 
$t_3$, $t_4$, $t_5$, $\dots$, with coefficients in 
the field $\mathbb{Q}(N)$ of rational functions 
in $N$. The size $N$ of matrices
is  considered as a variable here. The topology of this
formal power series ring is the Krull topology 
defined by setting
$\deg t_j = j$. The monomial
$\prod_j t_j ^{v_j(\Gamma)}$ is a finite product
for each connected
ribbon graph $\Gamma$, and has degree $2e(\Gamma)$. 
The matrix integral in LHS of Eqn.(\ref{eq:proto})
is meaningful only as an asymptotic series. 
Using the asymptotic technique developed in 
\cite{Mulase1995, Mulase1998}, we can obtain a
well-defined formal power series in 
$\mathbb{Q}(N)[[t_3,t_4,t_5,\dots]]$ 
from the
calculation of the integral. 
{F}irst we consider a truncation of variables
$(t_3,t_4,\dots,t_{2m})$ for some $m$. 
If we replace the integral
with
$$
Z(m) = \int_{\mathcal{H}_{N,\mathbb{C}}}
e^{-\frac{1}{2}N\tr X^2} 
e^{N\sum_{j= 3} ^{2m} \frac{t_j}{j}\tr X^j}d\mu(X)\;,
$$
then it converges and defines a holomorphic
function on 
$$
(t_3,\dots,t_{2m-1}, t_{2m})
\in \mathbb{C}^{2m-3}\times \{t_{2m}\in\mathbb{C}\;|
\; Re(t_{2m})<0\}\;.
$$
The holomorphic function has a unique asymptotic
series expansion at 
$$
(t_3,\dots,t_{2m-1}, t_{2m})=0\; .
$$
We can show that the  terms of degree
$n$ in the asymptotic
expansion of $Z(m)$ are stable for every $m> n$.
Therefore, $\lim_{m\rightarrow \infty} Z(m)$
is convergent in the Krull topology of 
$\mathbb{Q}(N)[[t_3,t_4,t_5,\dots]]$, and determines
a well-defined element. LHS of Eqn.(\ref{eq:proto})
is the logarithm of this limit.

Let us now consider a
finite-dimensional $*$-algebra $A$ together with 
a linear map called \emph{trace}
$$
\langle\; \rangle : A\longrightarrow \mathbb{C}
$$
satisfying that $\langle ab \rangle = \langle ba \rangle$.
A typical example is the group algebra $A = \mathbb{C}[G]$
of a finite group $G$. We define a $*$-operation by
$$
*: \mathbb{C}[G] \owns x=\sum_{\gamma\in G} x^\gamma 
\cdot \gamma \longmapsto x^* = 
\sum_{\gamma\in G}
\overline{x^\gamma} \cdot \gamma^{-1}\in \mathbb{C}[G]\;.
$$
As a trace, we use the character of the regular representation
$\rchi_{\reg}$, by linearly
extending it to the whole
group algebra. Let $\mathcal{H}_{\mathbb{C}[G]}$ 
denote the real vector subspace of $\mathbb{C}[G]$
consisting of self-adjoint elements. We need a 
Lebesgue measure on $\mathcal{H}_{\mathbb{C}[G]}$.

The self-adjoint condition $x^*=x$ means
$x^{\gamma^{-1}}=\overline{x^\gamma}$.
Let us decompose the group $G$ into the disjoint union
of three subsets
\begin{equation}
\label{eq:groupdecomp}
G=G_I \cup G_+ \cup G_-\; ,
\end{equation}
where $G_I = \{\gamma\in G\;|\;\gamma^2=1\}$ is the 
set of involutions of $G$, and $G_+$ and $G_-$ are
chosen so that
\begin{equation}
\label{eq:G+-}
\begin{cases}
G\setminus G_I = G_+\cup G_-\\
*(G_+) = G_-\; .
\end{cases}
\end{equation}
Every self-adjoint element of $\mathbb{C}[G]$ is written
as 
\begin{equation}
\label{eq:selfadjointC}
x = \sum_{\gamma\in G_I} x^\gamma \cdot \gamma 
+ \frac{1}{\sqrt{2}}
\sum_{\gamma\in G_+} 
\big(y^\gamma \cdot (\gamma + \gamma^{-1})
+ i z^\gamma \cdot (\gamma -\gamma^{-1})\big)\;,
\end{equation}
where $x^\gamma$, $y^\gamma$, and $z^\gamma$ are real numbers.
Thus we have 
$\mathcal{H}_{\mathbb{C}[G]} = \mathbb{R}^{|G|}$
as a real vector space. We see from Eqn.(\ref{eq:selfadjointC})
that if $x$ is self-adjoint, then
\begin{equation}
\label{eq:Ceucl}
\frac{1}{|G|}\rchi_\reg(x^2)
=\sum_{\gamma\in G_I}
(x^\gamma)^2 + \sum_{\gamma\in G_+} \big((y^\gamma)^2
+(z^\gamma)^2\big) \; .
\end{equation}
This is a non-degenerate quadratic form on $\mathbb{R}^{|G|}$.
Let $dx$ denote the Lebesgue measure on 
$\mathcal{H}_{\mathbb{C}[G]}$  that is invariant under the
Euclidean transformations with respect to the quadratic 
form (\ref{eq:Ceucl}).
A \emph{normalized} Lebesgue measure is defined by
\begin{equation}
\label{eq:normalizedC}
d\mu(x) = \frac{dx}{\int_{\mathcal{H}_{\mathbb{C}[G]}}
\exp\left(-\frac{1}{2}\rchi_\reg( x^2)\right) \;dx}\;.
\end{equation}

\begin{thm}
\label{thm:CGasymptotic}
As a formal power series in infinitely many variables
$t_3$, $t_4$, $t_5$, $\dots$, we have an equality
\begin{multline}
\label{eq:CGasymptotic}
\log \int_{\mathcal{H}_{\mathbb{C}[G]}}
\exp\left(-\frac{1}{2}\;\rchi_\reg(  x^2 )\right)
\exp\left(\sum_{j\ge 3} \frac{t_j}{j}\rchi_\reg(x^j )
\right) d\mu(x)\\
= \sum_{\substack{\Gamma  \text{ connected ribbon}\\
\text{graph with valence }\ge 3}}
\frac{1}{|\Aut_R \Gamma|}|G|^{\rchi(S_\Gamma)-1}
|\Hom(\pi_1(S_\Gamma),G)|
\prod_{j\ge 3} t_j ^{v_j(\Gamma)}\;,
\end{multline}
where $S_\Gamma$ is the oriented surface 
determined by a ribbon graph  $\Gamma$.
\end{thm}

\begin{proof}
In the computation of the integral, it is easier to use
a normalized trace function 
$\langle \;\rangle = \frac{1}{|G|} \rchi_\reg$
instead of $\rchi_\reg$. 
The formula to be established is then
\begin{multline}
\label{eq:CGasymptoticnormal}
\log \int_{\mathcal{H}_{\mathbb{C}[G]}}
\exp\left(-\frac{1}{2}\;\langle x^2 \rangle \right)
\exp\left(\sum_{j\ge 3} \frac{t_j}{j}\langle x^j \rangle 
\right) d\mu(x)\\
= \sum_{\substack{\Gamma  \text{ connected ribbon}\\
\text{graph with valence }\ge 3}}
\frac{1}{|\Aut_R \Gamma|}|G|^{f(\Gamma)-1}
|\Hom(\pi_1(S_\Gamma),G)|
\prod_{j\ge 3} t_j ^{v_j(\Gamma)}\;,
\end{multline}
where $f(\Gamma)$ is the number of faces in the 
surface $S_\Gamma$. 

{F}irst we expand the exponential factor:
$$
\exp\left(\sum_{j\ge 3} \frac{t_j}{j}\langle x^j \rangle
\right) = 
\sum_{v_3,v_4,v_5,\dots} \prod_{j\ge 3} 
\left(\frac{t_j ^{v_j}}
{v_j ! j^{v_j}} \langle x^j\rangle^{v_j}\right)\;.
$$
The coefficient of $\prod_j t_j ^{v_j}$ in the asymptotic
expansion is then 
\begin{equation}
\label{eq:asymptoticcoeff}
\prod_{j\ge 3} \frac{1}
{v_j ! j^{v_j}} 
\int_{\mathcal{H}_{\mathbb{C}[G]}}
\exp\left(-\frac{1}{2}\langle x^2 \rangle\right)
\prod_{j\ge 3} \langle x^j\rangle^{v_j} d\mu(x).
\end{equation}
We note that  (\ref{eq:asymptoticcoeff}) is 
a convergent integral. 
Let $y=\sum_{\gamma\in G} y^\gamma \cdot 
\gamma\in\mathbb{C}[G]$ be a variable running on  the
group algebra, and write 
$y^\gamma = u^\gamma + i w^\gamma$ as the sum of
real and imaginary parts. Now introduce
\begin{equation}
\label{eq:complexdy}
\frac{\partial}{\partial y} = \sum_{\gamma\in G}
\frac{\partial}{\partial y^\gamma}\cdot \gamma^{-1},
\end{equation}
where 
$$
\frac{\partial}{\partial y^\gamma} = 
\frac{1}{2}\;\left(\frac{\partial}{\partial u^\gamma}
- i\frac{\partial}{\partial w^\gamma}\right)\;.
$$

\begin{lem}
{F}or every $j >0$ and $n\ge 0$, we have
\begin{equation}
\label{eq:complexxy}
\left\langle \left(\frac{\partial}{\partial y}
\right)^j\right\rangle^n
e^{\langle x(y+y^*)\rangle} 
= \langle x^j\rangle^n \cdot 
e^{\langle x(y+y^*)\rangle} \;.
\end{equation}
\end{lem}

\begin{proof}
Since the $y$ derivative of $y^*$ is zero, we can ignore $y^*$ in 
the formula. 
\begin{equation*}
\begin{split}
\left\langle \left(\frac{\partial}{\partial y}
\right)^j\right\rangle
e^{\langle xy\rangle} &= 
\sum_{\gamma_1,\dots,\gamma_j}
\frac{\partial}{\partial y^{\gamma_1}}\cdots
\frac{\partial}{\partial y^{\gamma_j}}
\langle \gamma_1 ^{-1}\cdots \gamma_j ^{-1}\rangle 
e^{\sum_\gamma x^{\gamma^{-1}}y^\gamma}\\
&=\sum_{\gamma_1,\dots,\gamma_j}
x^{\gamma_1 ^{-1}}\cdots x^{\gamma_j ^{-1}}
\langle \gamma_1 ^{-1}\cdots \gamma_j ^{-1}\rangle 
e^{\langle xy\rangle} \\
&= \langle x^j\rangle e^{\langle xy\rangle}.
\end{split}
\end{equation*}
Applying this computation $n$ times, we obtain 
(\ref{eq:complexxy}).
\end{proof}
Using the completion of square, we have
\begin{equation}
\label{eq:comlexsquare}
-\frac{1}{2}\; \langle x^2\rangle + 
\langle x(y+y^*)\rangle
= -\frac{1}{2}\; \langle (x-(y+y^*))^2\rangle + 
\frac{1}{2}
\; \langle (y+y^*)^2\rangle.
\end{equation}
Thus
\begin{equation}
\label{eq:pair}
\begin{split}
\int_{\mathcal{H}_{\mathbb{C}[G]}}
\exp\left(-\frac{1}{2}\langle x^2 \rangle\right)
\prod_{j\ge 3} \langle x^j\rangle^{n} d\mu(x) 
&= \left.
\left\langle \left(\frac{\partial}{\partial y}\right)^j
\right\rangle^n
\int_{\mathcal{H}_{\mathbb{C}[G]}}
e^{-\frac{1}{2}\langle x^2\rangle + 
\langle x(y+y^*)\rangle }
d\mu(x)\right|_{y=0}\\
&= \left.
\left\langle \left(\frac{\partial}{\partial y}
\right)^j\right\rangle^n
\int_{\mathcal{H}_{\mathbb{C}[G]}}
e^{ -\frac{1}{2} \langle (x-(y+y^*))^2\rangle + 
\frac{1}{2}
 \langle (y+y^*)^2\rangle }
d\mu(x)\right|_{y=0}\\
&=\left.
\left\langle \left(\frac{\partial}{\partial y}
\right)^j\right\rangle^n
e^{ \frac{1}{2}
\langle (y+y^*)^2\rangle }\right|_{y=0} \;,
\end{split}
\end{equation}
where we used the fact that 
$y+y^*\in\mathcal{H}_{\mathbb{C}[G]}$
and the translational invariance of the Lebesgue measure
$d\mu(x)$.

\begin{lem}
We have
\begin{equation}
\label{eq:complexyder}
\frac{\partial}{\partial y} \; e^{ \frac{1}{2}
\langle (y+y^*)^2\rangle }
= (y+y^*) e^{ \frac{1}{2}
\langle (y+y^*)^2\rangle }\;.
\end{equation}
\end{lem}

\begin{proof}
\begin{equation*}
\begin{split}
\frac{\partial}{\partial y} \; e^{ \frac{1}{2}
 \langle (y+y^*)^2\rangle }
&= \sum_\gamma \frac{\partial}{\partial y^\gamma} \cdot
\gamma^{-1} \exp\left(\frac{1}{2}\sum_\gamma
(y^\gamma + \overline{y^{\gamma^{-1}}})
(y^{\gamma^{-1}} + \overline{y^\gamma})\right)\\
&=\sum_\gamma (y^{\gamma^{-1}} + \overline{y^\gamma})
{\gamma^{-1}}e^{ \frac{1}{2}
 \langle (y+y^*)^2\rangle }\\
&=(y+y^*)\; e^{ \frac{1}{2}
 \langle (y+y^*)^2\rangle } \;.
\end{split}
\end{equation*}
\end{proof}

The application of powers of $\langle (\partial/\partial y)^j\rangle$
to $e^{ \frac{1}{2}
\langle (y+y^*)^2 \rangle}$ produces zero contribution
unless two differentiations are paired
because of the restriction $y=0$ at the end. Therefore, 
(\ref{eq:pair})
counts the number of pairs of differential operators
in the product. Since the trace $\langle\; \rangle$ is invariant
under cyclic permutations, let us represent each
$\langle (\partial/\partial y)^j\rangle$ as a $j$-valent
vertex of a ribbon graph. Two vertices are connected
if the differentiations from the vertices are paired one another.
Since $y^*$ is killed by the differentiation and
\begin{equation}
\label{eq:edge}
\langle \cdots \frac{\partial}{\partial y}\cdots\rangle \; 
\langle \cdots y \cdots \rangle = 
\sum_{\gamma\in G}
\langle \cdots \gamma \cdots\rangle \; 
\langle \cdots \gamma^{-1} \cdots \rangle\;,
\end{equation}
the edge connecting the vertices
 comes with an assignment of a group element 
$\gamma$ and $\gamma^{-1}$ on the two half-edges.
With a factor $1/\prod_j v_j ! j^{v_j}$, the 
integral (\ref{eq:asymptoticcoeff}) is equal to
$$
\sum_{\substack{
\Gamma \text{ ribbon graph}\\
v_j(\Gamma) = v_j}}
\frac{1}{|\Aut_R \Gamma|} \mu_\Gamma(G) 
\prod_j t_j ^{v_j}\;,
$$
where $\mu_\Gamma(G)$ is the number of
 assignments of group elements to 
each half-edge of a ribbon graph $\Gamma$ subject to
the following conditions:

\medskip

\noindent
\emph{Condition} 1. If half-edges $E_+$ and $E_-$
form an edge $E$ of $\Gamma$ and a group element 
$w$ is assigned to $E_+$, then $w^{-1}$ is assigned to
$E_-$;

\smallskip

\noindent
\emph{Condition} 2. At every vertex, 
the product of all group elements
assigned to half-edges incident to the vertex according to 
the cyclic order of the vertex  is equal to $1$. 

\medskip

\noindent
The second condition comes from the fact
$$
\langle w_1 w_2\cdots w_j \rangle = 
\begin{cases}
1\qquad w_1 w_2\cdots w_j = 1\\
0\qquad \text{oherwise} \;,
\end{cases}
$$
that appears as a $j$-valent vertex of (\ref{eq:edge}).

\begin{lem}
\label{lem:topologicalinv}
The quantity $\mu_\Gamma(G)$ is a
topological invariant of a compact oriented surface
with a fixed number of marked points (or the faces
of its cell-decomposition).
\end{lem}

\begin{proof}
This follows from the invariance of $\mu_\Gamma(G)$
under an edge contraction and edge insertion \cite{MW}. 
When an edge that is not a loop is contracted, 
the configuration of group elements on the new graph  still
satisfies Conditions 1 and 2. If the edge is inserted back,
then we know exactly what group element has to be assigned
to each half-edge, due to Condition 2. 
This proves the lemma.
\end{proof}

\begin{figure}[htb]
\centerline{\epsfig{file=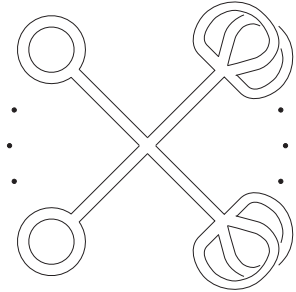, width=1.5in}}
\caption{A standard graph for a
closed  oriented surface of genus
$g$ with $f$ marked points, or faces. It has $f-1$ tadpoles
on the left and $g$ bi-petal flowers on the right.}
\label{fig:standard}
\end{figure}

As in \cite{MW}, we can use a standard graph 
for each topology to calculate the
number $\mu_\Gamma(G)$. If we use 
{F}igure~\ref{fig:standard}
as our standard graph $\Gamma$, 
then we immediately see that the number of
configurations of group elements 
on this graph is
$$
|G|^{f(\Gamma)-1}|\Hom(\pi_1(S_\Gamma), G)|\; ,
$$
where the tadpoles of the graph have the contribution
of $|G|^{f(\Gamma)-1}$ in the computation.

This completes the proof of the expansion formula
(\ref{eq:CGasymptoticnormal}), and by adjusting the
constant factor in the exponential function of the integrand,
 we establish
Theorem~\ref{thm:CGasymptotic}. 
\end{proof}

Note that we have a $*$-algebra isomorphism
\begin{equation}
\label{eq:complexdecomp}
\mathbb{C}[G] \isom \bigoplus_{\lambda\in\hat{G}}
\End(V_\lambda)\;,
\end{equation}
which decomposes the character of the regular representation
into the sum of irreducible characters: 
$$
\rchi_{\reg} = \sum_{\lambda\in\hat{G}} 
(\dim V_\lambda) \rchi_\lambda 
= \sum_{\lambda\in\hat{G}} N_\lambda \;
\tr_{V_\lambda}\;,
$$
where $N_\lambda = \dim V_\lambda$ and 
$\rchi_\lambda$ is its character.
Therefore, 
\begin{equation}
\label{eq:oricomputation}
\begin{split}
&\log \int_{\mathcal{H}_{\mathbb{C}[G]}}
\exp\left(-\frac{1}{2}\;\rchi_\reg(  x^2 )\right)
\exp\left(\sum_{j\ge 3} \frac{t_j}{j}\rchi_\reg (x^j )
\right) d\mu(x)\\
=&\log \int_{\mathcal{H}_{\mathbb{C}[G]}}
\prod_{\lambda\in\hat{G}}
\exp\left(-\frac{N_\lambda}{2}\;\tr_{V_\lambda}(  x^2 )\right)
\exp\left(N_\lambda \sum_{j\ge 3} 
\frac{t_j}{j}\tr_{V_\lambda}(x^j )
\right) d\mu_\lambda(x)\\
=&\sum_{\lambda\in\hat{G}}\log 
\int_{\mathcal{H}_{N_\lambda, \mathbb{C}}}
\exp\left(-\frac{N_\lambda}{2}\;\tr_{V_\lambda}(  x^2 )\right)
\exp\left(N_\lambda \sum_{j\ge 3} 
\frac{t_j}{j}\tr_{V_\lambda}(x^j )
\right) d\mu_\lambda(x)\\
=&\sum_{\substack{\Gamma  \text{ connected ribbon}\\
\text{graph with valence }\ge 3}}
\frac{1}{|\Aut_R \Gamma|}
\sum_{\lambda\in\hat{G}}
 (\dim V_\lambda)^{\rchi (S_\Gamma)}
\prod_{j\ge 3} t_j ^{v_j(\Gamma)}\;,
\end{split}
\end{equation}
where $d\mu_\lambda$ is the normalized Lebesgue
measure on the space of
$N_\lambda\times N_\lambda$ hermitian matrices.
Comparing the two expressions (\ref{eq:CGasymptotic})
and (\ref{eq:oricomputation}),
we obtain Eqn.(\ref{eq:orihomcount}).

\section{Non-orientable Case}

\noindent
{F}or the non-orientable case, 
the prototype integral is a matrix integral over $N\times N$
real symmetric
matrices \cite{BIPZ, GHJ}:
\begin{equation}
\label{eq:realproto}
\log  \int_{\mathcal{H}_{N,\mathbb{R}}}
e^{-\frac{1}{4}N\tr X^2} 
e^{\frac{N}{2}\sum_j \frac{t_j}{j}\tr X^j}d\mu(X)
=
\sum_{\substack{\Gamma \text{ connected}\\
\text{M\"obius graph}}}
\frac{1}{|\Aut \Gamma|} N^{\rchi(S_\Gamma)} 
\prod_{j} t_j ^{v_j (\Gamma)}\;,
\end{equation}
where what we call a \emph{M\"obius graph} is a  graph drawn on
a closed surface, orientable or non-orientable, defining
a cell-decomposition of the surface. Its automorphism is an
automorphism of the cell-decomposition of the surface that is
determined by the graph, but this time we allow 
orientation-reversing 
automorphisms.

Every compact non-orientable
surface without boundary
is obtained by removing $k$ disjoint disks 
from a sphere $S^2$ and glue $k$ \emph{cross-caps}
back into the holes. The number of cross-caps is called
the \emph{cross-cap genus} of the non-orientable surface.
If $S_\Gamma$ is non-orientable, then its cross-cap genus 
$k$ is determined by
$$
\rchi(S_\Gamma) = 2-k = v(\Gamma) -e(\Gamma) + 
f(\Gamma)\;,
$$
where again by $f(\Gamma)$ we denote the number of disjoint
open disks in $S_\Gamma\setminus \Gamma$.

The space $\mathcal{H}_{N,\mathbb{R}}$
of $N\times N$ real symmetric matrices is a real vector
space of dimension $N(N+1)/2$, and $d\mu(X)$ is the normalized
Lebesgue measure of this space. We note that the coefficients
of the integral in (\ref{eq:realproto}) are different from 
(\ref{eq:proto}), reflecting the fact that a dihedral group naturally 
acts on a vertex of a M\"obius graph.

We can generalize the matrix integral
(\ref{eq:realproto})  to an integral
over the real group algebra $\mathbb{R}[G]$, which is a
$*$-algebra with $\rchi_{\reg}$ as a trace function.

\begin{thm}
\begin{multline}
\label{eq:realgr}
\log  \int_{\mathcal{H}_{\mathbb{R}[G]}}
e^{-\frac{1}{4}\rchi_{\reg}( x^2)} 
e^{\frac{1}{2}\sum_j \frac{t_j}{j}
\rchi_{\reg}( x^j)}d\mu(x)\\
=
\sum_{\substack{\Gamma \text{ connected}\\
\text{M\"obius graph}}}
\frac{1}{|\Aut \Gamma|} |G|^{\rchi(S_\Gamma)-1} 
|\Hom(\pi_1(S_\Gamma),G)|
\prod_{j} t_j ^{v_j (\Gamma)}\;,
\end{multline}
where the integral is taken over the 
space of self-adjoint elements of
$\mathbb{R}[G]$,
and $S_\Gamma$ is the orientable or non-orientable
surface determined by a M\"obius graph $\Gamma$.
\end{thm}

{F}or the proof of this expansion formula, we refer to
\cite{MY}. 
Recall that the real group algebra $\mathbb{R}[G]$
decomposes into simple factors
according to the three types of irreducible representations
(\ref{eq:3reps}). 
{F}irst we note that $\hat{G}_1$ consists of complex irreducible
representations of $G$ that are defined over $\mathbb{R}$.
A representation in $\hat{G}_2$ is not defined over
$\mathbb{R}$, and its character is not real-valued. Thus
the complex conjugation acts on the set $\hat{G}_2$
without fixed points. Let $\hat{G}_{2+}$ denote a half
of $\hat{G}_2$ such that
\begin{equation}
\label{eq:G2+}
\hat{G}_{2+}\cup \overline{\hat{G}_{2+}} = \hat{G}_2\;.
\end{equation}
{F}inally, a complex irreducible 
representation of $G$ that belongs to $\hat{G}_3$
admits a skew-symmetric 
bilinear form. In particular, its dimension (over $\mathbb{C}$)
is even. Now we have a $*$-algebra isomorphism
\begin{equation}
\label{eq:realdecomp}
\mathbb{R}[G] \isom \bigoplus_{\lambda\in \hat{G}_1}
\End_{\mathbb{R}}(V_\lambda ^\mathbb{R})
\dsum \bigoplus_{\lambda\in \hat{G}_{2+}}
\End_{\mathbb{C}}(V_\lambda)
\dsum \bigoplus_{\lambda\in \hat{G}_3}
\End_{\mathbb{H}}(V_\lambda ^\mathbb{H})\;,
\end{equation}
where $V_\lambda ^\mathbb{R}$ is a real irreducible
representation of $G$ that satisfies 
$V_\lambda = V_\lambda ^\mathbb{R}\tensor_\mathbb{R}
\mathbb{C}$. The space $V_\lambda ^\mathbb{H}$ is a 
$(\dim_\mathbb{C} V_\lambda)/2$-dimensional vector
space defined over $\mathbb{H}$ such that its image 
under the natural injection
\begin{equation}
\label{eq:HtoC}
\End_\mathbb{H} (V_\lambda ^\mathbb{H})
\longrightarrow \End_\mathbb{C} (V_\lambda)
\end{equation}
coincides with the image of 
$$
\rho_\lambda : \mathbb{R}[G]\longrightarrow
\End_\mathbb{C} (V_\lambda)\;,
$$
where $\rho_\lambda$ is the representation of $G$ corresponding
to $\lambda\in\hat{G}_3$. The injective algebra homomorphism
Eqn.(\ref{eq:HtoC}) is defined by
$$
1\longmapsto \begin{pmatrix}
1\\&1\end{pmatrix}\;, \quad
i\longmapsto \begin{pmatrix}
&1\\-1\end{pmatrix}\;, \quad
j\longmapsto \begin{pmatrix}
i\\&-i\end{pmatrix}\;, \quad
k\longmapsto \begin{pmatrix}
&-i\\-i\end{pmatrix}\;.
$$
The algebra isomorphism Eqn.(\ref{eq:realdecomp})
gives a formula for the character
of the regular representation on $\mathbb{R}[G]$:
\begin{equation}
\label{eq:realreg}
\rchi_{\reg} = \sum_{\lambda\in\hat{G}_1} (\dim_\mathbb{R} 
V_\lambda ^\mathbb{R}) \rchi_\lambda +
 \sum_{\lambda\in\hat{G}_{2+}} (\dim_\mathbb{C} 
V_\lambda) (\rchi_\lambda +  \overline{\rchi_\lambda})+
\sum_{\lambda\in\hat{G}_3}2 (\dim_\mathbb{C} 
V_\lambda) \cdot \trace_{V_\lambda ^\mathbb{H}}\; ,
\end{equation}
where in the last term the character 
is given as the trace of quaternionic
$(\dim_\mathbb{C} V_\lambda)/2\times 
(\dim_\mathbb{C} V_\lambda)/2$ matrices.
To carry out the non-commutative integral
Eqn.(\ref{eq:realgr}), we need to know the
result for quaternionic self-adjoint matrix integrals. {F}ortunately, a
recent paper \cite{MW} provides exactly this necessary formula:
\begin{equation}
\label{eq:quaternionintegral}
\log \int_{\mathcal{H}_{N,\mathbb{H}}}
e^{ N \tr X^2} e^{2N \sum_j \frac{t_j}{j}\tr X^j d\mu(X)}
= \sum_{\substack{\Gamma \text{ connected}\\
\text{M\"obius graph}}}
\frac{1}{|\Aut \Gamma|}(-2N)^{\rchi(S_\Gamma)}
\prod_j t_j^{v_j(\Gamma)}\;,
\end{equation}
where $\mathcal{H}_{N,\mathbb{H}}$ is the space of 
$N\times N$ quaternionic self-adjoint matrices.
We pay a particular attention to the negative sign in RHS of the 
formula and the factor $2N$. The extra factor $2$ cancels
with the size of the quaternionic matrices
in $\End_\mathbb{H}(V_\lambda ^\mathbb{H})$,
which is half of the dimension of $V_\lambda$. 
The negative sign in RHS of (\ref{eq:quaternionintegral})
is the source of the negative sign in the second summation 
term of Eqn.(\ref{eq:nonorihomcount}).

The computation of the matrix integral of each factor
of Eqn.(\ref{eq:realdecomp})
then establishes 
Eqn.(\ref{eq:nonorihomcount}).
Note that the $\hat{G}_2$ component has no contribution in
this formula. This is due to the fact that graphical
expansion of a complex
hermitian matrix integral contains only oriented ribbon graphs.

\section{Algebraic Proof of the Counting Formula}

\noindent
The counting formulas 
(\ref{eq:orihomcount}) and (\ref{eq:nonorihomcount})
have an
 algebraic proof \cite{FS, Med} without going through the
computation of matrix integrals. Although the proof 
is well known to experts in group theory \cite{Jones},
since it is not found in modern textbooks, 
we record it here
to illuminate its relation to
our  graphical expansion formulas. The necessary
backgrounds for the algebraic proof are found in 
Isaacs \cite{Isaacs}, Serre \cite{Serre} and 
Stanley \cite{Stanley}.

Let us first identify the group algebra 
$$
\mathbb{C}[G] = \bigg\{x=\sum_{\gamma\in G} x(\gamma)\cdot
\gamma \bigg\}
$$
with the vector space $F(G)$ of functions on $G$. The
\emph{convolution product} of  two functions
$x(\gamma)$ and $y(\gamma)$ is defined by
$$
(x*y)(w) = \sum_{\gamma\in G}x(w\gamma^{-1})y(\gamma)\;,
$$
which makes $(F(G),*)$ an algebra 
 isomorphic to the group algebra. In this isomorphism,
the set of class functions $C\!F(G)$ corresponds to the center 
$Z\mathbb{C}[G]$ of 
$\mathbb{C}[G]$. 
According to the decomposition into simple factors 
Eqn.(\ref{eq:complexdecomp}), we have an algebra 
isomorphism
$$
Z\mathbb{C}[G] = \bigoplus_{\lambda\in\hat{G}}
\mathbb{C}\;,
$$
where each factor $\mathbb{C}$ is the center of
$\End{V_\lambda}$. The projection to each factor is
given by an algebra homomorphism 
$$
pr_\lambda : Z\mathbb{C}[G]\owns
x = \sum_{\gamma\in G} x(\gamma) \cdot \gamma
\longmapsto 
pr_\lambda(x)=\frac{1}{\dim V_\lambda}
\sum_{\gamma\in G} x(\gamma) \rchi_\lambda(\gamma)
\in\mathbb{C}\;.
$$
Now let
\begin{equation}
\label{eq:proj}
p_\lambda = \frac{\dim V_\lambda}{|G|}
\sum_{\gamma\in G} \rchi_\lambda (\gamma^{-1})\cdot
\gamma \in Z\mathbb{C}[G]\;.
\end{equation}
The orthogonality of the irreducible characters shows that
$pr_\lambda(p_\mu) = \delta_{\lambda\mu}$. 
Consequently, we have
$p_\lambda p_\mu = \delta_{\lambda\mu} p_\lambda$, or 
more concretely, 
\begin{equation*}
\begin{split}
&\frac{\dim V_\lambda}{|G|}
\sum_{s\in G} \rchi_\lambda (s^{-1})\cdot
s \cdot 
\frac{\dim V_\mu}{|G|}
\sum_{t\in G} \rchi_\mu (t^{-1})\cdot t\\
&=\frac{\dim V_\lambda\cdot \dim V_\mu}{|G|^2}
\sum_{w\in G}\left(
\sum_{t\in G} \rchi_\lambda ((wt^{-1})^{-1})
\rchi_\mu (t^{-1})\right)
\cdot w\\
&=\delta_{\lambda\mu}
\frac{\dim V_\lambda}{|G|}
\sum_{w\in G} \rchi_\lambda (w^{-1})\cdot w\;.
\end{split}
\end{equation*}
Therefore, we have
$$
\sum_{t\in G} \rchi_\lambda (tw^{-1})
\rchi_\mu (t^{-1})
=  \frac{|G|}{\dim V_\mu}\delta_{\lambda\mu}
\rchi_\lambda (w^{-1})\;,
$$
or in terms of the convolution product,
\begin{equation}
\label{eq:charconv}
\rchi_\lambda *\rchi_\mu = 
\frac{|G|}{\dim V_\mu}\delta_{\lambda\mu}
\rchi_\lambda\;.
\end{equation}

Now
let us consider
\begin{equation}
\label{eq:combinatorics}
\begin{cases}
f_g (w) =\big|
\{(a_1,b_1,a_2,b_2,\dots,a_g,b_g)
\in G^{2g}\;|
\; a_1b_1a_1 ^{-1}b_1 ^{-1}
\cdots a_g b_g a_g ^{-1} b_g ^{-1} =
w\}\big|\\
r_k (w) =\big|\{(a_1,a_2,\dots,a_k)\in G^{k}\;|
\; a_1 ^2 a_2 ^2\cdots a_k ^2 =
w\}\big|\;.
\end{cases}
\end{equation}
Since these are class functions
on the group $G$,  they can be written as  linear combinations of 
the irreducible characters 
$\{\rchi_\lambda\;|\; \lambda\in\hat{G}\}$
of $G$.

\begin{thm}
\label{thm:768+769}
{F}or every $g\ge 1$ and $w\in G$, we have
\begin{equation}
\label{eq:Stanley768}
f_g(w) = \sum_{\lambda\in \hat{G}} \left(
\frac{|G|}{\dim_\mathbb{C} V_\lambda}\right)^{2g-1}
\cdot \rchi_\lambda (w)\;.
\end{equation}
Corresponding to the non-orientable case, for every $k\ge 1$,
we have
\begin{equation}
\label{eq:Stanley769}
r_k(w) = \sum_{\lambda\in\hat{G}_1} \left(
\frac{|G|}{\dim_\mathbb{C} V_\lambda}\right)^{k-1}\cdot
\rchi_\lambda(w) 
-\sum_{\lambda\in\hat{G}_3} \left(-\;
\frac{|G|}{\dim_\mathbb{C} V_\lambda}\right)^{k-1}\cdot
\rchi_\lambda(w)\;. 
\end{equation}
\end{thm}

\begin{proof}
{F}rom the definition (\ref{eq:combinatorics}),
we see that
$$
f_{g_1 + g_2} = f_{g_1} * f_{g_2}\qquad
\text{and}\qquad
r_{k_1+k_2} = r_{k_1} * r_{k_2}\;.
$$
In particular, 
$$
f_g = \overset{g\text{-times}}{\overbrace{f_1 *\cdots * f_1}}
\qquad \text{and}\qquad
r_k = \overset{k\text{-times}}{\overbrace{r_1 *\cdots * r_1}}\;.
$$
The general formulas now follow from the formulas
for $f_1$ and $r_1$ that are found in \cite{Stanley}
and \cite{Isaacs}, respectively,
using (\ref{eq:charconv}).
\end{proof}

Evaluating Eqns.(\ref{eq:Stanley768}) and (\ref{eq:Stanley769})
at $w=1$, we obtain the counting formulas
(\ref{eq:orihomcount}) and (\ref{eq:nonorihomcount}).
The quantity $g$ in the formula is of course the genus of 
an oriented surface, and $k$ is the cross-cap genus of
a non-orientable surface. The class function $f_g$ (resp.\ 
$r_k$) gives the number $|\Hom(\pi_1(S),G)|$ 
when evaluated
at $w=1$ for an oriented (resp.\ non-orientable) surface $S$.

\section{Remarks}

\noindent
Hermitian matrix integrals have been used effectively
in the study of topology of moduli spaces of
Riemann surfaces with marked points 
\cite{Harer-Zagier, Kontsevich, 
 Penner, W1991a}. These works rely on a relation between
graphs and Riemann surfaces (cf.~\cite{Strebel}).
There is a burst of developments
in this direction lately \cite{O, OP1, OP2}. 
Real symmetric matrix integrals are used for
the study of  moduli spaces of real algebraic curves
\cite{GHJ}.  
In a context of finite groups, a striking relation between random
matrices and representation theory of symmetric
groups is discovered in \cite{BDJ, BDJ2, BOO, BR, D, J, O}. 
We have seen in this article that an extension of
these matrix integrals shows yet another interesting
relation between surface geometry and finite group theory.

There is a completely different  proof of the
 counting formula Theorem~\ref{thm:orientable}
due to Freed and Quinn \cite{FQ}.  They use
Chern-Simons gauge theory with a finite gauge group. To establish 
the formula, they 
construct a Chern-Simons gauge theory on each Riemann surface
$S$. Interestingly, our approach does not start with 
a specific manifold, which is usually the space-time in physics.
In a sense our non-commutative integral is a 
quantum field theory on a finite group \emph{without}
space-time, and surfaces
appear in its Feynman diagram expansion.
 Yet the result shows that Chern-Simons gauge theory
on a surface is mathematically equivalent to our
non-commutative integrals over group algebras. 
To be more precise, the group algebra model is a
generating function of Chern-Simons gauge theory with finite
gauge group for \emph{all} closed surfaces. 
The original algebraic proofs of the formulas using the
 convolution product found in \cite{FS,Med} are 
indeed related to the cut-and-paste construction 
of topological quantum field theory that reduces the
construction of an invariant on a surface of genus $g$ to
surfaces of lower genera \cite{FQ, W1991}. 

A more recent development on counting
formulas in finite groups is found in \cite{L3}.

\bigskip

\begin{ack}
The authors thank Andrei Okounkov for drawing their
attention to many classical literatures of the subject. 
They also thank Dmitry Fuchs,
Greg Kuperberg, Anne Schilling, Albert Schwarz,
Bill Thurston and Andrew Waldron for
many valuable comments and stimulating discussions
on the subjects of this paper. 
\end{ack}

\bigskip

%Bibliography

\providecommand{\bysame}{\leavevmode\hbox to3em{\hrulefill}\thinspace}

\bibliographystyle{amsplain}

\end{document}